\documentclass[preprint,12pt]{elsarticle}

\usepackage{lipsum}  
\usepackage{amsmath,amssymb,amsthm} 
\usepackage{graphicx}
\usepackage{bm}
\usepackage[normalem]{ulem}

\newcommand{\AAA}{\lambda} 
\newcommand{\BBB}{\mu} 
\newcommand{\CCC}{\nu} 
\newcommand{\MM}{\mathcal{M}} 
\newcommand{\ww}{\omega} 
\newcommand{\RR}{\mathbb{R}}
\newcommand{\NN}{\mathbb{N}}
\newcommand{\FORMATGAMMA}{\bm{\gamma}}
\newcommand{\POINTYALPHA}{\Lambda}

\journal{Nonlinear Waves: Computation and Theory-XI}

\begin{document}

\begin{frontmatter}



\title{Hasimoto Transformation of General Flows Expressed in the Frenet Frame}

\author[label1]{Jacob S.~Hofer}
\affiliation[label1]{organization={Department of Applied Mathematics and Statistics, Colorado School of Mines},
            addressline={1500 Illinois St.},
            city={Golden},
            postcode={80401-1887},
            state={CO},
            country={USA}
           }
\ead{jhofer@mines.edu}
\author[label1]{Scott A.~Strong}
\ead{sstrong@mines.edu}



\begin{abstract}
A one-dimensional space curve in $\mathbb{R}^{3}$ is a useful nonlinear medium for modeling vortex filaments and biological soft-matter capable of supporting a variety of wave motions. 
Hasimoto’s transformation defines a mapping between the kinematic evolution of a space curve and nonlinear scalar equations evolving its intrinsic curve geometry. This mapping is quite robust and able to transform general vector fields expressed in the Frenet frame, resulting in a fully nonlinear integro-differential evolution equation, whose coefficient structure is defined by the coordinates of the flow in the Frenet frame.  In this paper, we generalize the Hasimoto map to arbitrary flows defined on space curves, which we test against several existing kinematic flows. After this, we consider the time dynamics of length and bending energy to see that binormal flows are generally length preserving, and bending energy is fragile and unlikely to be conserved in the general case. 
\end{abstract}




\end{frontmatter}


\section{Introduction}
The autonomous motions of slender filaments are essential to several areas of applied science, e.g., atmospheric, aerodynamic, oceanographic phenomena, astrophysical plasmas, and superfluid and superconducting states of matter~\cite{Andersen2014IntroductionEquilibrium}. Through the Frenet-Serret equations, the filament's one-dimensional centerline configuration is prescribed by two scalar functions, i.e., curvature and torsion. The curvature and torsion dynamics determine the evolution of a filament, which can be found by evolving the Frenet-Serret coordinate system itself. However, an auxiliary condition is required to close the system of equations. Often the kinematics of the filament are known and act as the constitutive relation connecting our geometric response to a physical context. Consequently, this framework seeks to reduce the problem of understanding filament behavior through its curvature and torsion dynamics. To aid in this, we have Hasimoto's remarkable result~\cite{Hasimoto1972}, which complexifies the local curvilinear coordinates so that the kinematics can be married to the Frenet-Serret differential equations resulting in a complex-valued partial differential equation controlling the evolution of the centerline geometry.

Hasimoto's transformation, initially established in the context of fluid mechanics where vortex lines served as an object for modeling and simulation, maps a global prescription for the motion of the filament centerline to a complex-valued scalar evolution of its curvature and torsion and continues to see widespread use with 38\% of its citations occurring in the last ten years. This is partly due to the applicability of vortex filaments in the setting of quantum liquids and Bose-Einstein condensation, the latter being a hotbed of activity since their experimental inception 25 years ago~\cite{Levi2001CornellCondensates, PNAS2014QT}. Here the filaments act as the skeleton of quantum turbulence~\cite{Tsubota2013QuantizedTurbulence} but they need not be exotic structures. Currently, it is hypothesized that they provide a natural object for singularity formation in Navier-Stokes \cite{Brenner2016PotentialEquations} and mediate the anomalous dissipation conjectured by Onsager~\cite{Eyink2006OnsagerTurbulence}, with their helical wave motion stitching together our picture of turbulent cascades~\cite{Walmsley2014DynamicsSpectrab}. 

While these techniques provide a dimensional reduction to fluid problems, i.e., fluid flow reduced to the evolution of a complex scalar representation of a space curve, they are fundamentally geometric in character and when coupled to the elastic mechanics pioneered by Euler and Bernoulli, similar analyses yield possible geometric configurations of stiff polymers structuring DNA supercoiling, self-assembly of bacterial fibers, actin filaments, and cell flagella~\cite{Ludu2012NonlinearSurfaces}. In either fluid or elastic contexts, the primary purpose of Hasimoto's transformation is to recast the behavior specified by a  setting's kinematics to the nonlinear evolution of intrinsically geometric quantities  where wave motion is less opaque. 

In the following, we generalize the Hasimoto transform to centerline evolutions expressed in the Frenet-Serret basis with arclength and time dependent coefficients. We show that the transformations of known special cases, e.g., local induction approximation~\cite{Hama1988GenesisLIA,Ricca1996TheDynamicsb}, generalized local induction equation~\cite{Strong2012GeneralizedTurbulence, Fukumoto1991Three-dimensionalVelocity}, and stiff chain polymers~\cite{Ludu2012NonlinearSurfaces}, are special cases of our general scalar evolution law. Additionally, we show specific conditions for which arclength is conserved and that conservation of bending energy is fragile. The work is organized into three sections. First, we outline the procedure necessary to arrive at the main result, which is a nonlinear integro-differential equation in evolutionary form that specifies the curve implicitly through its curvature and torsion dynamics. Next, we apply this result to known kinematic flows of space curves. After this, we consider the relationship between the kinematic flows expressed in the Frenet-Serret basis and the global quantities of arclength and bending energy, before we conclude.

\section{Generalization of Hasimoto Transformation to Frenet Framed Kinematics}

We consider the Hasimoto transformation of a one-dimensional space curve $\FORMATGAMMA : \RR^{1+1} \rightarrow \RR^3$, parameterized by arclength, $s$, and permitted to evolve in time, $t$, where a background velocity field defines the kinematic equation 
\begin{align}
    \FORMATGAMMA_t &=  \CCC \,  \textbf{T} + \BBB \, \textbf{N} + \AAA \, \textbf{B} \label{flow law}
\end{align}
expressed in the local curvilinear basis defined by the tangent, $\textbf{T}$, normal $\textbf{N}$, and binormal, $\textbf{B}$, vectors whose coordinates $\AAA$, $\BBB$, and $\CCC$, are arbitrary functions of arclength and time. This coordinate system implicitly defines a space curve via its curvature, $k$, and torsion, $\tau$, describing local rotations about the binormal and tangent vectors, respectively. Hasimoto's transformation begins by defining  a change of coordinates $\MM = (\textbf{N}+i\textbf{B})  e^{i\theta}$ where $\theta = \int^s \tau  ds'$, which is used in conjunction with a complex-valued wave function whose modulus defines curvature and derivative of phase defines torsion, i.e., $\psi = k e^{i\theta}$.

Through the Frenet-Serret equations, $\textbf{T}_s =k\textbf{N}$, $\textbf{N}_s = -k\textbf{T} + \tau \textbf{B}$, and $\textbf{B}_s = -\tau \textbf{N}$, it is possible to derive the following expressions for the derivatives of $\MM$ and $\textbf{T}$ with respect to arclength,
\begin{align}
    \MM_s &= - \psi \textbf{T} \label{Ms}\\
    \textbf{T}_s &= \frac{1}{2} (\overline{\psi} \MM + \psi \overline{\MM}),
\end{align}
where we have denoted complex conjugation as $\overline{\psi} = k e^{-i\theta}$. Furthermore, we have that the time dynamics of the Hasimoto frame are given by the partial derivative of $\MM$ with respect to time, 
\begin{align}
    \MM_t = iR\MM + \ww \textbf{T} \label{Mt}
\end{align}
where $R$ is an arbitrary real function and $\omega$ is defined by projection with the standard(real) inner-product, 
\begin{align}
    \ww &= \langle \MM_t,\textbf{T}\rangle = -\langle\MM,\textbf{T}_t\rangle = -(\BBB_s + i \tau \BBB + i ( \AAA_s + i \tau \AAA) + \CCC k) e^{i \theta} \\&=-i \left( \frac{\AAA}{k}\psi \right) _s - \left(\frac{\BBB}{k}\psi\right)_s -\CCC \psi,
\end{align}
since $\langle \MM , \textbf{T} \rangle_t =0$. Equating the mixed partial derivatives of $\MM$ maps the curve dynamics onto the Hasimoto wave function and requires a representation for the partial derivative of our kinematic equation, Eq.~(\ref{flow law}), 
\begin{align}
    \textbf{T}_t &= (\AAA_s + \BBB \tau ) \textbf{B} + (\BBB_s - \AAA \tau + \CCC k) \textbf{N} + (\CCC_s - \BBB k)\textbf{T}  \label{t_t},
\end{align}
which can be expressed in the orthogonal basis of $\MM$, $\overline{\MM}$, and $\textbf{T}$ as, 
\begin{align}
    \textbf{T}_t = -\frac{1}{2}(\ww \overline{\MM} + \overline{\ww} \MM)  + (\CCC_s - \BBB k) \textbf{T}.
\end{align}
Consequently, we can compute the mixed partials of Eq.~(\ref{Ms}) and Eq.~(\ref{Mt}) to find,
\begin{align}
    \MM_{ts} &= i R_s \MM + \frac{\ww}{2}\left(\overline{\psi} \MM + \psi \overline{\MM}\right) - i R \psi \textbf{T}\\
    \MM_{st} &= \frac{\psi}{2}(\overline{\ww} \MM + \ww \overline{\MM}) -\psi_t \textbf{T} + \left(\CCC_s - \BBB k\right) \textbf{T} 
\end{align}
and equating the coefficients of the $\textbf{T}$ and $i \MM$ terms gives
\begin{gather}
    R_s = \frac{1}{2}i(\ww \overline{\psi} - \overline{\ww} \psi) \label{Rs} \\
    \psi_t - i R \psi + \ww_s + \psi(\CCC_s - \BBB k) =0 \label{that}.
\end{gather}
Integrating Eq.~(\ref{Rs}) with respect to arclength yields,
\begin{align}
    R= \AAA k - \POINTYALPHA + \int^s {\BBB \tau k }\, ds' \label{R},
\end{align}
where $\POINTYALPHA_k = \AAA$ and 
%
combining Eq.~(\ref{that}) and Eq.~(\ref{R}) results in a fully nonlinear integro-differential equation 
\begin{align}
    i \psi_t + \left( \frac{\AAA}{k}\psi  - i \frac{\BBB}{k}\psi\right)_{ss} - i \CCC \psi_s + \left(\AAA k - i\BBB k  - \POINTYALPHA + \int^s {\BBB \tau k }\, ds' \right) \psi =0 \label{pde}
\end{align}
While complicated, an equation of this form is not surprising since perturbations of simpler kinematic laws have revealed similar consequences~\cite{Majda2002VorticityFlow}. It is interesting to note that while the initial and boundary conditions are defined by the curve geometry, i.e., curvature and torsion, the coefficients defining the evolution of $\psi$ are found by specifying a flow, Eq.~(\ref{flow law}). 

\section{Scalar Nonlinear Evolution Equations for Various Kinematic Forms}

Originally, Hasimoto used this procedure to transform a curvature dependent binormal flow, known as the vortex filament equation (VFE), or local induction approximation/equation,  
\begin{align}
    \FORMATGAMMA_t = k \textbf{B}, \label{binormal flow}
\end{align}
which corresponds to our generalized flow where $\AAA = k$, $\BBB = 0$, and $\CCC =0$. Thus, simplifying Eq.~(\ref{flow law}) and expressing our coefficients in terms of $k=|\psi|$,  Eq.~(\ref{pde}) reduces to 
%
\begin{align}
    i \psi_t + \psi \left(k^2 - \frac{k^2}{2}  \right) +  \left( \psi \right) _{ss} =
    i \psi_t + \frac{1}{2}|\psi|^2 \psi +   \psi  _{ss}  =0
\end{align}
which is an integrable cubic focusing nonlinear Schr\"odinger equation, supporting soliton transport~\cite{Hasimoto1972}, breathing modes~\cite{Salman2013BreathersVortices}, and nonlinear dispersion~\cite{Newton1987StabilityWaves}, implying that even in the simplest cases, vortex filaments in ideal fluids are a nonlinear medium. 

The VFE, Eq.~(\ref{binormal flow}), defines a specific flow arising from an asymptotic analysis of the Biot-Savart integral~\cite{Lamb1980ElementsTheory} and is the simplest example of a class of kinematic equations controlling the evolution of curve centerlines. In 1991, Fukumoto and Miyazaki~\cite{Fukumoto1991Three-dimensionalVelocity} developed a correction to VFE accounting for axial velocity, 
\begin{align}
    \FORMATGAMMA_t = k \textbf{B} + W\left(\frac{1}{2} k^2 \textbf{T} + k_s \textbf{N} +k \tau \textbf{B}\right)
\end{align}
where $W$ is a constant. This corresponds to Eq.~(\ref{flow law}) with $\AAA = k + W k \tau$, $\BBB = W k_s$, and $\CCC = W \frac{1}{2} k^2$. Using these coefficients we find that Eq.~(\ref{pde}) can be reduced to
\begin{align}
    i \psi_t + \frac{1}{2}|\psi|^2 \psi + \psi_{ss}  -i W\left(  \psi_{sss} + \frac{3}{2} |\psi|^2 \psi_s\right) =0,
\end{align}
which is an integrable mixture of the previous nonlinear Schr\"odinger equation and the modified Korteweg-De Vries equation showing, again, that vortex filaments are fundamentally connected to nonlinear wave motion. 

While the integro-differential equation, Eq.~(\ref{pde}) is a general result, it is impractical to analyze. At the same time, we have several results where specific flows reduce to  differential equations.  The following provides an example of how the integral terms can be circumnavigated by the introduction of mild auxiliary assumptions, e.g., series expansions in the geometric variables, so that nonlinear partial differential equations result. For example, if we generalize binormal flow to be a nonlinear function of the curvature variable, 
\begin{align}
    \FORMATGAMMA_t = \alpha(k) \textbf{B},
\end{align}
which is supported by a more detailed asymptotic analysis of the Biot-Savart integral~\cite{Strong2012GeneralizedTurbulence},
 then $\AAA = \alpha(k)$, and $\BBB= \CCC =0$ and 
%
\begin{align}
    i \psi_t +   \left( \frac{\AAA}{k}\psi \right) _{ss} + \psi \left( \AAA k - \POINTYALPHA \right)  =0. \label{strongintegro}
\end{align}
However, since $\POINTYALPHA = \int \alpha dk$, Eq.~(\ref{strongintegro}) is still in intregro-differential form.
If we assume that the magnitude of the binormal flow can be expanded in powers of curvature,
%
\begin{align}
    \alpha(k) = \sum_{n=1}^{\infty} a_n k^n,
\end{align}
where $a_{n} \in \mathbb{R}$ for $n=1,2,3,\dots$, then we can explicitly calculate the necessary anti-derivative 
\begin{align}
    \POINTYALPHA = \int \alpha\, dk= \sum_{n=1}^\infty \frac{a_n}{n+1} k^{n+1}.
\end{align}
Substituting these series into Eq.~(\ref{strongintegro}) gives
%
\begin{align}
    i \psi_t +   \left( \frac{\AAA}{k}\psi \right) _{ss} + f(|\psi|) \psi  =0,
\end{align}
where 
\begin{align}
    f(|\psi|) =  \sum_{n=1}^{\infty} a_n \frac{n}{n+1} |\psi|^{n+1},
\end{align}
 which is a fully nonlinear partial differential equation of Schr\"odinger type. Numerical investigations have shown that the nonlinear curvature terms endow the medium with enhanced dispersion coupled with a gain/loss mechanism not seen in VFE and appears to support the redistribution of bending energy \cite{Strong2017Non-HamiltonianCondensates}. 


In the previous examples, we have shown that our Hasimoto transformation of a generalized kinematic equation governing the flow of $\FORMATGAMMA$ reproduces known results from the fluid dynamics literature. That said, space curves are applicable models in other settings. In particular, through Euler's elastica theory one can consider the nonlinear dynamics of stiff polymers, which is applicable to a variety of biological settings. Indeed, by specifying the forces acting on the polymer in the Frenet-Serret basis and minimizing an elastic potential energy functional, one can arrive at a set of equations describing the nonlinear dynamics about a static equilibrium shape~\cite{Ludu2012NonlinearSurfaces}. Applying the Hasimoto transformation to these equations yields an integro-differential equation consistent with Eq.~(\ref{pde}). Indeed, formal connections between rigid body motions and ideal fluid dynamics as geodesic flows over Riemannian manifolds exist~\cite{Arnold1966SurParfaits}. Consequently, these geometric connections through the Hasimoto map are not unexpected, and worth continued investigation. 

\section{Characterization of Conserved Quantities Under Example Kinematic Prescriptions}
With our generalization of the Hasimoto transformation, we can now relate the kinematic flow to  two quantities important in the characterization of space curve dynamics, i.e., the total length and bending energy. Considering the dynamics of these two quantities, we derive conditions for which they are conserved by the autonomous flow affecting the curve geometry, Eq.~(\ref{flow law}). First, we integrate over the arclength domain, $D=[0,L]$, to define the length of $\bm{\gamma}$ as, 
\begin{align}
    I_1 = \int_D \left\lVert \FORMATGAMMA_s \right\rVert  \, ds \label{arclength}. 
\end{align}
Under the arclength parameterization, the time rate of change of this quantity is given by,  
\begin{align}
    \frac{dI_1}{dt} &= \int_{0}^L \frac{\partial}{\partial t} \left\lVert \FORMATGAMMA_s \right\rVert  \, ds \\
    &= \int_{0}^L  (\FORMATGAMMA_{st} \cdot \FORMATGAMMA_s) \, ds\\
    & = \int_{0}^L ( \textbf{T}_{t} \cdot \textbf{T} ) \, ds \\ 
    & = \int_{0}^L (\CCC_{s} - \BBB k )\, ds.\label{arceq2}
\end{align}
%

If the flow is purely binormal, i.e., $\CCC=0$ and $\BBB=0$, then this integral is zero, and thus arclength is an invariant of the flow. Such fields result from the local induction approximation, which are known to correspond to metric preserving Killing vector fields~\cite{Langer1991PoissonEquation}.  While the transformation of our curvature dependent binormal flow complicates the evolution of the curve geometry, the arclength is conserved by virtue of binormality of the flow. Importantly such fields are consistent with the Hamiltonian description of ideal fluid dynamics which treats the Euler equations as arclength conserving Hamiltonian dynamic ~\cite{Khesin2012SymplecticMembranes}. 

More generally, we find that conservation is maintained if $\CCC_s - \BBB k=0$.  This is particularly interesting because it defines conditions where flow in the tangential and normal directions can interact in order to conserve total length, giving us insight into those flows consistent with the Eulerian theory. While the flow law of Fukumoto and Miyazaki is far from simplistic, it  too  satisfies this condition implying that it is also Killing. It is unclear whether this is a consequence of their reductions to the Biot-Savart integrand~\cite{Fukumoto1991Three-dimensionalVelocity}, which defines an approximation that avoids the use of elliptic representations, or is generally true in filament settings. 
The case of bending energy is notably more difficult.

There is a direct relationship between energy stored in a space curve and the square of its state variables, in this case, curvature.  The total bending energy along the curve is given by, 
\begin{align}
    I_2 = \int_D k^2 \, ds \label{benergy}
\end{align}
whose time derivatives is 
\begin{align}
    \frac{d I_2}{dt} &= \int_0^L 2 k_t k \, ds.
\end{align}
The dynamics of curvature are implicit as the real part of Eq.~(\ref{pde}). However, it is easiest to consider the time derivative of the Frenet-Serret equations coupled to Eq.~(\ref{flow law}). From here, we can project out the curvature component to get
\begin{align}
    \frac{d I_2}{dt} &= 2\int_0^L  \left[  k \BBB_{ss} + 2k^2 \CCC_s + k k_s \CCC -2k \AAA_s \tau - k\AAA\tau_s - k \BBB \tau^2 -  \BBB k^2 \tau \right] ds. 
\end{align}

For the VFE we have that when $\AAA = k$ and $\BBB=\CCC =0$, and the integrand is the exact derivative of $-k^2 \tau$, which vanishes for periodic curves. Generally, however, the dynamics of bending energy are highly dependent on the form of the flow, and unlikely to be conserved.
%
For example, letting $\AAA = k^n$ for some $n \in \NN$ when $\nu=\mu=0$, and keeping in mind periodic boundary conditions, the integral becomes 
\begin{align}
    \frac{d I_2}{dt} &= \int_0^L -4k (k^n)_s \tau - 2k(k^n)\tau_s ds \\
    &= \int_0^L \frac{2n-2}{n+1} {k^{n+1}} \tau_s ds
\end{align}
which we can see vanishes when $n=1$, i.e., the case of VFE. Generally, however, conservation must come from a delicate interplay between $k$ and $\tau$, which seems to be violated in simulations~\cite{Strong2018Non-HamiltonianCondensates}.



\section{Conclusion and Outlook}
As the forces responsible for the autonomous evolution of space curves are often found through approximation, one should expect a diverse assortment of possible centerline kinematics. In this paper, we have generalized Hasimoto's transformation for flows represented in the Frenet-Serret basis. The resulting fully nonlinear integro-differential equation is not unexpected and though it is cumbersome, we find that important special cases yield descriptions of the autonomous dynamics as evolutionary partial differential equations. Additionally, its generality allows us to connect the kinematic prescription of curve evolution to geometrically conserved quantities and integrability. Arclength is robust for binormal flows. That said, the loss of bending energy conservation allows for a coupling between the ambient space responsible for the kinematics and the filament itself. This insight may help bring resolution to problems related to phonon emission in Bose-Einstein condensates and its relationship with the Kelvin wave cascade in quantum turbulence~\cite{PhysRevB.90.094501}.

\section*{Acknowledgements}

The authors would like to thank the Colorado School of Mines Summer Undergraduate Research Fellowship for supporting this project.




 \bibliographystyle{model1b-num-names.bst} 
 \bibliography{2018Bib}


\end{document}